\baselineskip=15pt plus 2pt
\magnification =1200
\def\sqr#1#2{{\vcenter{\vbox{\hrule height.#2pt\hbox{\vrule width.#2pt
height#1pt\kern#1pt \vrule width.#2pt}\hrule height.#2pt}}}}
\def\square{\mathchoice\sqr64\sqr64\sqr{2.1}3\sqr{1.5}3}
 at 10truept
\font\smalletters=cmr8 at 10truept
 at 10truept
 at 10truept at 12truept

\font\medtenrm=cmr10 scaled\magstep2
{\smalletters JMAA-07-1648 revised}\par
\centerline {\medtenrm Integrable operators and the squares of Hankel
operators}\par
\vskip.05in
\centerline {\medtenrm Gordon Blower}\par
\vskip.05in
\centerline  {\sl Department of Mathematics and Statistics, 
Lancaster University}\par
\centerline  {\sl Lancaster, LA1 4YF, England, UK. E-mail: 
g.blower@lancaster.ac.uk}\par
\vskip.05in
\centerline {5th September 2007}\par
\vskip.05in
\vskip.05in
\hrule
\vskip.05in
{{\noindent {\bf Abstract} Integrable operators arise in random matrix theory,
where they describe the asymptotic eigenvalue distribution of large self-adjoint random
matrices from the generalized unitary ensembles. This paper gives sufficient conditions
for an integrable operator to be the square of a Hankel operator, and applies the
condition to the Airy, associated Laguerre, modified Bessel and Whittaker functions.\par
\vskip.05in
\noindent {\sl Keywords:} random matrices, Tracy--Widom operators
\vskip.05in
\noindent MSC2000 Classification 47B35\par
\vskip.05in
\vskip.05in
\hrule
\vskip.05in
\vskip.05in
\noindent {\bf 1. Introduction}\par
\vskip.05in
\indent Integrable operators with kernels of the form 
$$W(x,y)={{f(x)g(y)-f(y)g(x)}\over{x-y}}\eqno(1.1)$$
have applications in quantum field theory
and random matrix theory, where they are used to describe the asymptotic distribution
of large random matrices; see [5, 17, 18]. Tracy and Widom [19] observed that many 
important distributions in random
matrix theory can be defined using solutions of systems
$$m(x){{d}\over{dx}}\left[\matrix{f(x)\cr g(x)\cr}\right] =
\left[\matrix{\alpha (x)&\beta (x)\cr -\gamma (x)&
-\alpha (x)\cr}\right]\left[\matrix{f(x)\cr g(x)\cr}\right]\eqno(1.2)$$
\noindent where $m(x), \alpha (x),\beta (x)$ and $\gamma (x)$ are polynomials. In
[17, 18], Tracy and Widom considered the Airy and Bessel kernels which
describe the soft and hard edges of generalized unitary ensembles, and 
proved the apparently miraculous identities that
the operators with these kernels were squares of self-adjoint Hankel operators; then they used this property to
compute their eigenfunctions and eigenvalues. \par
\indent Let $K$ be a separable Hilbert space, and 
$L^2((0, \infty );dx; K)$ be the Bochner--Lebesgue space of strongly
measurable functions $\phi: (0, \infty )\rightarrow K$ such that $\int_0^\infty \Vert \phi (x)\Vert^2_Kdx<\infty$.

\vskip.05in
\noindent ------------------\par
\vskip.05in
\noindent {\smalletters This work was partially supported by EU Network Grant MRTN-CT-2004-511953
`Phenomena in High Dimensions'.}\par
\vfill
\eject
\noindent In [2] we considered a general class
of differential equations which gives rise to integrable operators that have
the form $W=\Gamma^*\Gamma$, where 
$\Gamma :L^2(0, \infty )\rightarrow L^2((0, \infty ); K)$ is a
continuous Hankel operator. The general
theorem of [2] specialised to the Airy and Bessel kernels, and in this paper 
we prove related results which deal with other integrable operators.\par
\indent After recalling some definitions, we state and prove the main Theorem
1.2, and then in section 2 give some applications. In section 3 we discuss the scope of Theorem 1.2 as it applies to (1.2) in reduced forms, and then in section 4 prove another theorem
which encompasses other applications, as in section 5.\par
\vskip.05in
\noindent {\bf Definition} {\sl (Integrable operators)} Let $I$ be a subinterval of
${\bf R}$. An integrable operator on $L^2(I; dx)$ is a continuous linear operator $W$ with kernel 
$$W(x,y)=2\sum_{j=1}^n {{f_j(x)g_j(y)}\over{x-y}}\qquad (x,y\in I; x\neq y)\eqno(1.3)$$
\noindent where $f_j,g_j$ are bounded and measurable functions such that $\sum_{j=1}^n
f_j(x)g_j(x)=0$ almost everywhere on $I$.\par
\vskip.05in
\noindent {\bf Lemma 1.1.} {\sl Suppose further that the $f_j$ and $g_j$ are
real-valued. Then $W$ is self-adjoint if and only if 
$$W(x,y)={{\bigl\langle Jv(x), v(y)\rangle}\over{x-y}}\qquad (x,y\in I; x\neq
y)\eqno(1.4)$$
\noindent where} $v(x)={\hbox{col}}[f_1(x), \dots ,f_n(x); g_1(x), \dots , g_n(x)]$
{\sl and} 
$$J=\left[\matrix{0& -I_n\cr I_n&0\cr}\right]\eqno(1.5)$$
\noindent {\sl with identity matrix $I_n\in M_n({\bf R})$ and the usual inner
product on ${\bf R}^{2n}$.}\par
\vskip.05in
\noindent {\bf Proof.} Clearly the kernel of $W$ is symmetric if and only if the numerator of 
$W(x,y)$ is skew-symmetric, in which case we can write 
$$2\sum_{j=1}^n f_j(x)g_j(y)=\sum_{j=1}^n {f_j(x)g_j(y)}-\sum_{j=1}^n
{f_j(y)g_j(x)},\eqno(1.6)$$
\noindent and the matrix expression follows directly.\par
\rightline{$\square$}\par
\vskip.05in 
\noindent {\bf Definition} {\sl (Hankel operator)} Let $\phi\in L^2((0, \infty ); K)$. The Hankel operator with symbol
$\phi$ is the integral operator 
$$\Gamma_\phi f(s)=\int_0^\infty \phi (s+t)f(t)\, dt.\eqno(1.7)$$
\vskip.05in 
\noindent Nehari's theorem [13] gives a sufficient condition for 
$\Gamma_\phi :L^2(0,\infty)\rightarrow L^2((0, \infty ); K)$ to be
continuous and gives an expression
for the operator norm $\Vert \Gamma\Vert$; while the Hilbert--Schmidt norm of
$\Gamma_\phi$ satisfies 
$$\Vert\Gamma_\phi \Vert_{HS}^2=\int_0^\infty s\Vert \phi (s)\Vert^2_K\, ds\eqno(1.8)$$
\noindent when $\sqrt s\phi (s)\in L^2((0, \infty ); K).$ Clearly
$\Gamma_\phi$ is self-adjoint when $K={\bf R}$.\par
\vskip.05in 
\indent Given an integrable operator, it is often valuable to identify a Hankel
operator $\Gamma_\phi$ such that $\Gamma_\phi^*\Gamma_\phi =W$ and to determine
whether $W$ is of trace class. In particular, when
$K={\bf R}$ and $W=\Gamma_\phi^2$, the spectral resolution
of the self-adjoint operator $\Gamma_\phi$ determines the spectral resolution of $W$.
This is the basis of the successful calculations in [17, 18, 19], which also exploited the fact
the eigenvectors of $\Gamma_\phi$ can be relatively easy to analyze. In [12], 
Megretski\u\i , Peller and Treil characterized the spectral
multiplicity function of a self-adjoint Hankel operator. Further, in applications to
determinantal point fields as in [15], one often wishes to show that $W$ is of trace class and
satisfies $0\leq W\leq I$.\par
\vskip.05in
\indent We consider first the case associated with the differential equation
$${{dv}\over{dx}}=J(\Omega_1x+\Omega_0+\Omega_{-1}x^{-1})v.\eqno(1.9)$$
\vskip.05in
\noindent {\bf Theorem 1.2.} {\sl Suppose that $\Omega_1, \Omega_0$ and $\Omega_{-1}$ are
real symmetric $(2n)\times (2n)$ constant matrices such that $\Omega_1\geq 0$ and 
$-\Omega_{-1}\geq 0$. Suppose further that $v$ satisfies (1.9), and that $v(x)$ and
$v(x)/x$ are bounded functions in $L^2((0, \infty );dx; {\bf R}^{4n}).$ Then there
exists a real linear subspace $K$ of ${\bf R}^{4n}$ with} ${\hbox{dim}}(K)\leq 
{\hbox{rank}}(\Omega_1)+{\hbox{rank}}(\Omega_{-1})$ {\sl and 
$\phi\in L^2((0, \infty );dx;K)$ such that $\Gamma_\phi :L^2(0, \infty
)\rightarrow L^2((0, \infty );dx;K)$ is continuous
and the kernel 
$$W(x,y)= {{\langle Jv(x), v(y)\rangle}\over {x-y}}\eqno(1.10)$$ 
\noindent factors as}
$$W=\Gamma_\phi^*\Gamma_\phi.\eqno(1.11)$$
\indent {\sl In particular, if} ${\hbox{rank}}(\Omega_1)+{\hbox{rank}}(\Omega_{-1})=1$,
{\sl then}
$W=\Gamma_\psi^2$ for some $\psi\in L^2(0,\infty )$.} 
\vskip.05in
\noindent {\bf Proof.} By Lemma 2.1 of [2], which essentially depends upon the
continuity of the Hilbert transform on $L^2({\bf R})$, we know that $W$ is also a
continuous linear operator on $L^2(0, \infty )$.\par
\indent We have 
$$\Bigl({{\partial}\over{\partial x}} +{{\partial}\over{\partial y}}\Bigr)
W(x,y)={{1}\over {x-y}}\Bigl(\Bigl\langle J{{dv}\over {dx}}, v(y)\Bigr\rangle 
+\Bigl\langle J v(x),{{dv}\over {dy}}\Bigr\rangle \Bigr)\qquad (x\neq y),\eqno(1.12)$$
\noindent where $dv/dx$ and $dv/dy$ satisfy (1.9). Now $J^2=-I$ and $J^*=-J$, so the 
matrices involved in the differential equation satisfy
$$J^2(\Omega_1x+\Omega_0+\Omega_{-1}x^{-1})+(\Omega_1y+\Omega_0+\Omega_{-1}y^{-1})J^*J$$
$$=-\Omega_1x+ \Omega_1y-\Omega_{-1}x^{-1}+\Omega_{-1}y^{-1};\eqno(1.13)$$
\noindent so by dividing by $x-y$, we obtain 
$$\Bigl({{\partial}\over{\partial x}} +{{\partial}\over{\partial y}}\Bigr)
W(x,y)=-\bigl\langle \Omega_1 v(x), v(y)\bigr\rangle +\Bigl\langle\Omega_{-1} {{v(x)}\over{x}}, 
{{v(y)}\over{y}}\Bigr\rangle .\eqno(1.14)$$
\noindent We introduce the positive roots of the positive semidefinite matrices
$\Omega_1$ and $\Omega_{-1}$, and the column vector 
$\phi (x)={\hbox{col}}[\sqrt{\Omega_1}v(x),\sqrt{-\Omega}_{-1}v(x)/x];$
\noindent so that, $\phi\in L^2((0, \infty ); dx;{\bf R}^{4n})$ and 
$$\Bigl({{\partial}\over{\partial x}} +{{\partial}\over{\partial y}}\Bigr)
W(x,y)=-\langle \phi (x), \phi (y)\rangle \qquad (x,y >0).\eqno(1.15)$$
\noindent Integrating this equation, we obtain
$$W(x,y)=\int_0^\infty \langle \phi (x+t), \phi (y+t)\rangle\, dt+g(x-y)\eqno(1.16)$$
\noindent for some differentiable function $g$; but $W(x, y)$ and the integral converge
to $0$ as $x\rightarrow\infty $ or $y\rightarrow\infty$; so $g=0.$ Hence 
$W=\Gamma_\phi^*\Gamma_\phi$, and we deduce that $\Gamma_\phi$ is a continuous 
Hankel operator from $L^2(0, \infty )\rightarrow L^2((0,\infty );dx; {\bf R}^{4n}).$ \par
\indent Finally, we observe that $\phi$ takes values in a linear subspace $K$ of\par
\noindent ${\hbox{range}}(\Omega_1)\oplus {\hbox{range}}(\Omega_{-1})$ which has dimension less
than or equal to ${\hbox{rank}}(\Omega_1)+{\hbox{rank}}(\Omega_{-1}).$\par
\indent If the sum of the ranks of $\Omega_1$ and $\Omega_{-1}$ equals one, then
$\phi$ takes values in a one-dimensional real linear subspace of ${\bf R}^{4n}$, so
$\phi (x)=e\psi (x)$ for some unit vector $e\in {\bf R}^{4n}$ and some $\psi\in L^2(0, \infty )$; hence
$\Gamma_\psi$ is self-adjoint and $W=\Gamma_\psi^2$.\par
\rightline{$\square$}\par

\vskip.05in
\noindent {\bf Remarks 1.3.} (i) In their applications to random matrix
theory, Tracy and Widom considered integrable operators on unions of
intervals such as $\cup_{j=1}^m [a_{2j-1}, a_{2j}].$ Many analytical
problems reduce to considering one interval at a time, and so are
addressed by the current paper. In a subsequent article [20], they
generalized their results to kernels of the form 
$${{\langle C\varphi (x), \varphi (y)\rangle}\over{x-y}}$$
\noindent where $C$ is an antisymmetric $n\times n$ matrix.\par  
\indent (ii) Theorem 1.2 has an analogue for discrete kernels on ${\bf Z}_+$,
 as in [11].\par 
\vskip.05in
\noindent {\bf 2. Factorization for some differential equations with simple poles}\par
\vskip.05in
\noindent In this section we consider how Theorem 1.2 applies to some differential
equations that are satisfied by familiar special functions.\par
\vskip.05in 

\noindent {\bf 2.1. The Airy equation} \par
\indent The Airy function ${\hbox{Ai}}$ gives rise to a solution
$u(x)= {\hbox{Ai}}(x+s)$ of the differential equation
$$u''(x)=(s+x)u(x)\qquad (x>0).\eqno(2.1)$$
\noindent Since the standard asymptotic formula for the Airy function [16, p. 18] gives
$${\hbox{Ai}}(x)={{1}\over {2x^{1/4}\sqrt{\pi}}}\Bigl( 1+O(x^{-3/2}\Bigr)\exp\Bigl(
-{{2}\over {3}}x^{3/2}\Bigr)\qquad (x\rightarrow\infty )\eqno(2.2)$$
\noindent the hypotheses of Theorem 1.2 are satisfied by $v(x)={\hbox{col}}[\phi (x),
\phi'(x)]$ where $\phi (s)$\par
\noindent $={\hbox{Ai}}(x+s);$
so $W=\Gamma_\phi^2$ where $\Gamma_\phi$ is a Hankel operator in the Hilbert--Schmidt
class. See [2, 17] for more details.\par
\vskip.05in
\noindent {\bf 2.2. The Laguerre equation}\par
\indent  The Laguerre equation [14] may be expressed as 
$$u''(x)+\Bigl( -{{1}\over{4}}+{{n+1}\over{x}}\Bigr) u(x)=0,\eqno(2.3)$$
\noindent with solution $u(x)=xe^{-x/2}L_n^{(1)}(x)$ where 
$$L_n^{(1)}(x)={{x^{-1}e^x}\over {n!}} {{d^n}\over{dx^n}}\bigl(
x^{n+1}e^{-x}\bigr)\qquad (x>0)
\eqno(2.4)$$
\noindent is the Laguerre polynomial of degree $n$ and parameter $\alpha =1$. 
The Laplace transform of $u$ is the rational function
$${\cal L}(u; \lambda )=(n+1){{(\lambda -{{1}\over{2}})^n}
\over {(\lambda +{{1}\over{2}})^{n+2}}}\qquad
(\Re\lambda >-1/2).$$
\noindent Theorem 1.2 applies directly to the system
$${{d}\over{dx}}\left[\matrix{u(x)\cr u'(x)\cr}\right]=
\left[\matrix{0&1\cr {{1}/{4}}-{{(n+1)}/{x}}&0\cr}\right]\left[\matrix{u(x)\cr u'(x)\cr}\right]
\eqno(2.5)$$ 
\noindent and gives the formula
$${{u'(x)u(y)-u'(x)u(y)}\over{x-y}}=(n+1)\int_0^\infty {{u(x+t)u(y+t)}\over
{(x+t)(y+t)}}\, dt\eqno(2.6)$$
\noindent where $\phi (x)=u(x)/x$ gives a Hankel operator $\Gamma_\phi$ of
Hilbert--Schmidt type.\par
\vskip.05in
\noindent {\bf 2.3. The Bessel equation}\par 
\indent The differential equation 
$$ u''(x)+{{1}\over{x}}u(x)=0\eqno(2.7)$$
\noindent has solution $u(x)=\sqrt {x}J_1(2\sqrt{x})$, where $J_1$ is the Bessel
function of the first kind of order one. The Laplace transform of $u$ satisfies ${\cal
L}(u, \lambda )=\lambda^{-2}\exp (-1/\lambda )$. Now the standard asymptotic formula for the
Bessel function [10, p. 171] shows that 
$$\phi (x)= {{u(x)}\over{x}}\asymp {{2^{1/4}}\over {\sqrt{\pi}}} x^{-3/4} \cos \Bigl(2\sqrt
{x}-{{3\pi}\over{4}}\Bigr)\qquad (x\rightarrow\infty );\eqno(2.8)$$
\noindent so $\phi$ belongs to $L^2(0, \infty )$; hence one can follow the proof of Theorem 1.2 and derive the formula
$${{u(x)u'(y)-u'(x)u(y)}\over {x-y}}=\int_0^\infty{{u(x+t)u(y+t)}\over
{(x+t)(y+t)}}\, dt\qquad (x,y>0).\eqno(2.9)$$
\noindent Here the Hankel operator $\Gamma_\phi$ is not
Hilbert--Schmidt.\par
\vskip.05in
\noindent {\bf 2.4. The Carleman operator with multiple spectrum}\par
\vskip.05in
\indent  The system
$${{d}\over{dx}}\left[ \matrix{f\cr g\cr}\right] = \left[ \matrix{0 &1/x\cr 0&0 \cr}
\right] \left[ \matrix{f\cr g\cr}\right]\eqno(2.10)$$
\noindent has the form considered in Theorem 1.2 and evidently has solution $f(x)=\log x$ and $g(x)=1$; further
$$W(x,y)={{\log x-\log y}\over{x-y}}=\int_0^\infty
{{dt}\over{(x+t)(y+t)}}\qquad(x,y>0)\eqno(2.11)$$
\noindent has a similar form to an integrable kernel, except that $\log x$ is
unbounded. Power showed that Carleman's
operator $\Gamma$, where
$$\Gamma h(x)=\int_0^\infty {{h(y)dy}\over {x+y}}\qquad (h\in L^2(0,
\infty )),\eqno(2.12)$$
\noindent is continuous on $L^2(0, \infty )$ and has spectrum $[0, \pi ]$ with 
spectral multiplicity two; see [13]. Hence $W$ is a
continuous linear operator on $L^2(0, \infty )$ with spectrum $[0, \pi^2]$ with 
multiplicity two. This example illustrates that simple differential equations can give 
positive definite Hankel operators with multiple spectra.\par
\vfill\eject
\noindent {\bf 2.5. Parabolic cylinder functions: non factorization}\par
\vskip.05in
\indent For $\Re p>-1$, let $D_p$ be the parabolic cylinder function, which satisfies
$$D_p''(x)+\Bigl(p+{{1}\over{2}}-{{x^2}\over{4}}\Bigr)D_p=0,$$
\noindent and let
$$H(x,y)={{D_p(x)D'_p(y)-D_p'(x)D_p(y)}\over {x-y}}.\eqno(2.13)$$
\noindent Then $\pm H$ is not the square of a self-adjoint Hankel operator
$\Gamma_\phi$. By
following the proof of Theorem 1.2, we obtain 
$$\Bigl({{\partial }\over {\partial x}}+{{\partial }\over {\partial
y}}\Bigr)H(x,y)=-{{1}\over{2}}(x+y)D_p(x)D_p(y),\eqno(2.14)$$
\noindent where $(x+y)D_p(x)D_p(y)/2$ cannot equal $\pm \phi
(x)\phi (y)$; indeed, for suitable $x_1, x_2>0$, the $2\times 2$ matrix
$[x_j+x_k]_{j,k=1,2}$ has both 
positive and negative eigenvalues. When $n$ is a nonnegative integer, $D_n$ is known as
a Hermite function, and may be written  
$$\phi_n(x)=(n!)^{-1/2}(2\pi )^{-1/4}(-1)^n e^{x^2/4} {{d^n}\over {dx^n}} 
e^{-x^2/2}.\eqno(2.15)$$
\noindent Aubrun [1] considers self-adjoint Hankel operators 
$\Gamma_{\phi_n}$ and $\Gamma_{\phi_{n+1}}$ such
that\par
\noindent ${H}=(1/2)(\Gamma_{\phi_n}\Gamma_{\phi_{n+1}}+
\Gamma_{\phi_{n+1}}\Gamma_{\phi_n})$; this gives information about
the singular numbers of $H$. For Hankel squares one has more
precise information about the eigenvalues. Borodin and Okounkov [4] have considered the discrete Hermite kernel, and
derived the formula
$${{\phi_{m+1}(s)\phi_{n}(s)-\phi_{m}(s)\phi_{n+1}(s)}\over{m-n}}=
\int_s^\infty \phi_m(t)\phi_n(t)\, dt\qquad (m,n=0, 1, \dots , m\neq n).\eqno(2.16)$$
\noindent Here the variable $n$ is the degree of the Hermite polynomial factor in 
$\phi_n$, and (2.16) is essentially different from (1.11).\par 
\vskip.05in

\noindent {\bf 3. Reducing to standard form}\par
\vskip.05in
\noindent {\bf Definition} {(\sl Operator monotone)} Let $I$ be an interval in ${\bf
R}$. A continuous function $\omega
:I\rightarrow {\bf R}$ is operator monotone increasing if, whenever $S$ and $T$ are
continuous and self-adjoint linear operators on Hilbert space that have spectra in $I$, 
$$S\leq T\Rightarrow  \omega (S)\leq \omega (T).\eqno(3.1)$$
\vskip.05in
\indent Further, $\omega$ is operator monotone if and only if the matrices 
$$\Bigl[{{\omega (x_j)-\omega (x_k)}\over {x_j-x_k}}\Bigr]_{j,k=1,
\dots , m},\eqno(3.2)$$
\noindent with diagonal entries $\omega' (x_j)$, are positive semidefinite for all $m=2, 3, \dots$ and $x_j\in I$ with $j=1,
\dots ,m$. In [2] we used
Loewner's characterization of operator monotone functions, which shows in particular
that an operator monotone function on $(0, \infty )$ extends to an analytic function
on a domain $U$ containing $(0, \infty )$ as in [9, p. 541].\par
\indent Theorem 1.2 shows that, under mild technical conditions, $W$ admits of a 
factorization $W=\Gamma_\phi^*\Gamma_\phi$ whenever 
$\omega (x)=\langle \Omega (x)\xi, \xi\rangle$ is operator monotone 
 on $(0, \infty )$, for all $\xi\in {\bf
R}^{4n}$  where $\Omega (x)=\Omega_1x+\Omega_0+\Omega_{-1}x^{-1}$.\par
\indent Suppose that $I$ has $0$ as an endpoint, and
let $U$ be a domain that contains $I.$
Suppose that $A(z)$ is a matrix function into $M_{2n}({\bf C})$ that is analytic on
$U$, except for an isolated singularity at $z=0$, and that 
$${{dv}\over{dz}}=A(z)v\quad (z\in U).\eqno(3.3)$$
\noindent By a standard change of variable, we mean $w(z)=T(z)v(z)$, where the analytic
function 
$T:U\rightarrow M_n({\bf C})$ has $T(z)$ invertible as a matrix for each $z\in U$.
The following result is commonly known as Birkhoff's normal form, although the
first correct statement and proof is due to Turrettin [21]. Gantmacher considered
some related examples which resemble (2.10) in [6, p. 146].\par
\vskip.05in
\noindent {\bf Proposition 3.1.} {\sl Suppose that
$A(z)=\sum_{k=-\infty}^{-1}A_kz^k$ is a Laurent expansion that converges for all $z\neq
0$. Then there exists a standard change of variable that reduces (3.3) to 
$${{dw}\over{dz}}=\Bigl( {{A_{-1}}\over {z}}+ {{A_{-2}}\over {z^2}}\Bigr) w.\eqno(3.4)$$
\indent Further, if the eigenvalues $\lambda_j$ have differences $\lambda_j-\lambda_k$
that are never equal to a natural integer, then one can remove the term in $A_{-2}$.}\par
\vskip.05in
\indent The appearance of the term $A_{-2}$ in (3.4) is important, since $-1/x^2$
is not operator monotone on $(0, \infty )$ by [9, p. 554]. So we cannot
simply adapt the proofs in section 1 to deal with the case in which $A_{-2}$ appears. 
However, if $\omega$ is operator monotone on $(0,
\infty )$, then $\omega (\sqrt{t})$ is likewise. This suggests the change of
independent variable $x=\sqrt{t}$, which we exploit in the examples in section 5.
Further, we adjust the definition of the kernel and the Hankel operators so that we
can obtain a factorization theorem in the next section.\par
\vskip.05in
\noindent {\bf Definition} {(\sl Hankel operator)} For $I=(1, \infty ),$ we use the Hankel operator 
$$\Gamma_\psi g(x)=\int_1^\infty\psi (xy)g(y){{dy}\over{y}}\eqno(3.5)$$
\noindent where $\psi\in L^2((1, \infty ); dy/y; K)$; whereas for $I=(0,1)$ we use
$$\Gamma_\rho h(x)=\int_0^1\rho (xy)h(y){{dy}\over{y}}\eqno(3.6)$$
\noindent where $\rho\in L^2((0,1); dy/y; K)$. These definitions reduce to the   
case $I=(0, \infty)$ in section 1 by the changes of variables $y=e^t$ and $y=e^{-t}$ respectively.\par
\vskip.05in

\noindent {\bf 4. Factorization theorem for differential equations with double poles}\par
\vskip.05in
\indent In this section we consider the differential equation 
$$x{{dv}\over{dx}}=J(\Omega_1x+\Omega_0+\alpha J+\Omega_{-1}x^{-1})v.\eqno(4.1)$$
\noindent To accommodate forthcoming examples, we have introduced the skew-symmetric 
matrix $\alpha J$ into
the constant term for some $\alpha\in {\bf R}$.\par
\vskip.05in   
\noindent {\bf Theorem 4.1.} 
{\sl Suppose that $\Omega_1, \Omega_0$ and $\Omega_{-1}$ are
real symmetric $(2n)\times (2n)$ constant matrices such that $\Omega_1\geq 0$ and 
$-\Omega_{-1}\geq 0$. Suppose further that $v$ satisfies (4.1), and that
$x^{\alpha }v(x)$ and
$x^{\alpha -1}v(x)$ are bounded functions in $L^2((1, \infty );dx; {\bf R}^{4n}).$ Then there
exists a real linear subspace $K$ of ${\bf R}^{4n}$ with} ${\hbox{dim}}(K)\leq 
{\hbox{rank}}(\Omega_1)+{\hbox{rank}}(\Omega_{-1})$ {\sl and $\phi\in L^2((1, \infty );
 dx/x;K)$ such that 
$\Gamma_\phi :L^2((1, \infty ); dx/x)\rightarrow L^2((1, \infty );dx/x ;K)$ 
is continuous and the kernel 
$$W(x,y)= {{(xy)^{(2\alpha +1)/2}}\over {x-y}}\bigl\langle Jv(x),
v(y)\bigr\rangle\eqno(4.2)$$
\noindent factors as}
$$W=\Gamma_\phi^*\Gamma_\phi.\eqno(4.3)$$
\noindent {\sl In particular, if} ${\hbox{rank}}(\Omega_1)+{\hbox{rank}}(\Omega_{-1})=1$,
{\sl then
$W=\Gamma_\psi^2$ for some $\psi\in L^2((1, \infty ); dx/x)$.}

\vskip.05in 
\noindent {\bf Proof.} We observe that by homogeneity
$$\Bigl( x{{\partial }\over {\partial x}}+ y{{\partial }\over {\partial y}}\Bigr)
{{(xy)^{(2\alpha +1)/2}}\over {x-y}}= 2\alpha {{(xy)^{(2\alpha +1)/2}}\over {x-y}}\qquad (x,y>0; x\neq y),\eqno(4.4)$$
\noindent and hence 
$$\eqalignno{\Bigl( x{{\partial }\over {\partial x}}+ y{{\partial }\over {\partial y}}\Bigr)W(x,y)
&= 2\alpha {{(xy)^{(2\alpha +1)/2}}\over {x-y}}\bigl\langle Jv(x), v(y)\bigr\rangle\cr
&\quad+
{{(xy)^{(2\alpha +1)/2}}\over {x-y}}\Bigl( \Bigl\langle Jx{{dv}\over{dx}}, v(y)\Bigr\rangle 
+\Bigl\langle Jv(x),y{{dv}\over{dy}}\Bigr\rangle\Bigr),&(4.5)}$$
\noindent  where the matrices involved in the final terms in (4.5) are
$$J^2(\Omega_1 x+\Omega_0+\alpha J+\Omega_{-1}x^{-1})+(
\Omega_1 y+\Omega_0-\alpha J+\Omega_{-1}y^{-1})J^*J$$
$$=-\Omega_1(x-y)+\Omega_{-1}(x-y)/(xy)-2\alpha J.$$
\noindent By cancelling the terms that involve $J$, we obtain
$$\Bigl( x{{\partial }\over {\partial x}}+ y{{\partial }\over {\partial
y}}\Bigr)W(x,y)=-(xy)^{2\alpha +1)/2}\langle \Omega_1v(x), v(y)\rangle+ (xy)^{(2\alpha +1)/2}\langle
\Omega_{-1}v(x), v(y)\rangle.\eqno(4.6)$$
\noindent We introduce the column vector
$$\psi (x)=\left[\matrix{\sqrt{\Omega_1}x^{(2\alpha +1)/2}v(x)\cr 
\sqrt{-\Omega}_{-1}x^{(2\alpha -1)/2}v(x)\cr}\right]\eqno(4.7)$$
\noindent which belongs to $L^2((1, \infty ); dx/x;{\bf R}^{4n})$ and satisfies
$$\Bigl( x{{\partial }\over {\partial x}}+ y{{\partial }\over {\partial
y}}\Bigr)\int_1^\infty \langle \psi (tx), \psi (ty)\rangle \, {{dt}\over{t}}=-\langle
\psi (x), \psi (y)\rangle .$$
\noindent Hence 
$$W(x,y)=\int_1^\infty\langle \psi (tx), \psi (ty)\rangle \, {{dt}\over{t}} +h(x/y)
\eqno(4.8)$$
\noindent where $h(x/y)\rightarrow 0$ as $x\rightarrow\infty$ or $y\rightarrow\infty$; so
$h=0.$ One can conclude the proof by arguing as in Theorem 1.2.\par
\rightline{$\square$}
\vskip.05in
\indent We consider later some examples in which $\Omega_{-1}=0$. In this case, we
can invoke the following existence theorem for solutions. \par
\vskip.05in

\noindent {\bf Proposition 4.2.} {\sl Suppose that the residue matrix $A_{-1}$ has eigenvalues $\lambda_j$
such that the differences $\lambda_j-\lambda_k$ are never equal to a natural integer.
Then the differential equation
$$z{{d}\over{dz}}X=(A_0z+A_{-1})X\eqno(4.9)$$
\noindent with $X(z)\in M_{2n}({\bf C})$ has a non-trivial solution of the form
$X(z)=Y(z)z^{A_{-1}}$, where $Y$ is an entire matrix function of order one.}\par
\vskip.05in

\indent For a proof see [8], where
Hille also discusses the asymptotic form of the solutions in terms of the Laplace
transform. Note that the Laplace transform of (4.9) has a similar form to (4.9)
itself and
in particular has the residue matrix $A_{-1}+I$.\par
\vskip.05in
\noindent {\bf 5. Examples of factorization for differential equations with double poles}\par
\vskip.05in
\noindent {\bf 5.1. Modified Bessel functions}\par
\vskip.05in
\indent For $0\leq \nu <1$, MacDonald's function is defined by
$$K_\nu (z)=\int_0^\infty e^{-z\cosh t}\cosh (\nu t)\, dt\qquad (\Re z>0)\eqno(5.1)$$ 
and satisfies the modified Bessel equation
$z^2K_\nu''+zK_\nu'-(\nu^2+z^2)K_\nu=0$; hence $u(x)=\sqrt {x}K_\nu (2\sqrt{x})$ satisfies 
$$u''(x)=\Bigl({{1}\over{x}}+{{\nu^2-1}\over{4x^2}}\Bigr)u(x).\eqno(5.2)$$ 
\noindent By [7, 8.451], $K_\nu (x)$ decays exponentially as $x\rightarrow\infty$. We can apply Theorem 4.1 to the system
$$x{{d}\over{dx}}\left[\matrix{u\cr w\cr}\right]=
\left[\matrix{-1&1\cr x-2+{{1}\over{4}}(\nu^2-1)&2\cr}\right]
\left[\matrix{u\cr w\cr}\right],\eqno(5.3)$$
\noindent so that, in terms of Theorem 4.1, 
$$\Omega_1=\left[\matrix{1&0\cr 0&0\cr}\right], \quad 
\Omega_0=\left[\matrix{-2+{{1}\over{4}}(\nu^2-1)&{{3}\over{2}}\cr {{3}\over{2}}&-1\cr}\right],\quad \alpha =-1/2, \quad
\Omega_{-1}=0,\eqno(5.4)$$
\noindent where the residue matrix $A_{-1}=J\Omega_0-\alpha I$ has eigenvalues
$(1/2)\pm (\nu /2)$.\par
\indent Thus one obtains
$${{u(x)v(y)-u(y)v(x)}\over {x-y}}=\int_1^\infty u(tx)u(ty) {{dt}\over{t}}$$
\noindent and after some reduction, one deduces that
$${{K_\nu (2\sqrt{x})\sqrt{y}K'_\nu (2\sqrt{y})-\sqrt{x}K'_\nu (2\sqrt{x}) 
K_\nu (2\sqrt{y})}\over {x-y}}=\int_1^\infty K_\nu (2\sqrt{tx}) 
K_\nu (2\sqrt{ty})dt\eqno(5.5)$$ 
\noindent where the right-hand side is the square of a Hankel operator of
Hilbert--Schmidt class.\par
\vskip.05in
\vskip.05in
\noindent {\bf 5.2. Bessel functions}\par
\vskip.05in
\indent The Bessel function $J_\nu$ satisfies
$x^2J_\nu''+xJ_\nu'+(x^2-\nu^2)J_\nu=0$, and hence $u=\sqrt {x}J_\nu (2\sqrt{x})$ satisfies 
$$u''(x)+\Bigl({{1}\over{x}}+{{1-\nu^2}\over{4x^2}}\Bigr)u(x)=0.\eqno(5.6)$$ 
\noindent One can apply Theorem 4.1 with some obvious sign changes to the system
$$x{{d}\over{dx}}\left[\matrix{u\cr w\cr}\right]=
\left[\matrix{-1&1\cr -x-2-{{1}\over{4}}(1-\nu^2)&2\cr}\right]
\left[\matrix{u\cr w\cr}\right],\eqno(5.7)$$
\noindent and after some reduction one obtains an identity from [18]
$${{\sqrt{x}J'_\nu (2\sqrt{x})J_\nu (2\sqrt{y})-J_\nu (2\sqrt{x}) 
\sqrt{y}J'_\nu (2\sqrt{y})}\over {x-y}}=\int_0^1 J_\nu (2\sqrt{tx}) 
J_\nu (2\sqrt{ty})dt.\eqno(5.8)$$ 
\vskip.05in
\vfill
\eject
\noindent {\bf 5.3. Whittaker's functions}\par
\vskip.05in
\noindent The homogeneous confluent hypergeometric equation may be reduced to
Whittaker's equation
$$w''+\Bigl( -{{1}\over{4}}+{{\kappa}\over {x}}+{{{{1}\over{4}}-\nu^2}\over {x^2}}
\Bigr)w=0,\eqno(5.9)$$
\noindent and the solutions of this are known as Whittaker's functions.
In [7, 9.227], the authors give a solution $w(z)=W_{\kappa , \nu }(z)$ such that
$w(z)\asymp e^{-z/2}z^\kappa$ as $z\rightarrow \infty$ along $(0, \infty )$. Kernels
involving Whittaker's function appear in [3].\par
\vskip.05in
\noindent {\bf Proposition 5.3.} {\sl Suppose that $u(x)=W_{\kappa , \nu }(2\sqrt{x})$
for some $\kappa \leq 0$. Then there exists a function $\Phi\in L^2((1, \infty
); dx /x; K)$ for a separable Hilbert space $K$ such that 
$$W(x,y)=(xy)^{1/4}{{u(x)u'(y)-u'(x)u(y)}\over {x-y}}\eqno(5.10)$$
\noindent factors as $W=\Gamma_\Phi^*\Gamma_\Phi$.}\par
\vskip.05in
\noindent {\bf Proof.} We can write, after a little reduction
$$x{{d}\over{dx}}\left[\matrix{u\cr v\cr}\right]={{1}\over{4}}\left[\matrix{ -4& \quad 4\cr 
x-2\kappa
\sqrt{x} -({{1}\over{4}}-\nu^2)-6&\quad 6\cr}\right]\left[\matrix{u\cr v\cr}\right],\eqno(5.11)$$
\noindent which gives 
$$\Omega (x)={{1}\over{4}}\left[\matrix{ x-2\kappa 
\sqrt{x} -({{1}\over{4}}-\nu^2)-6 &5\cr 5&-4}\right] -{{1}\over{4}}\left[\matrix{0&-1\cr
1& 0\cr}\right],\eqno(5.12)$$
\noindent hence $\alpha =-1/4$. The function $\sqrt x$ is operator monotone increasing 
on $(0, \infty )$; indeed, $\phi_t(x)=t^{1/4}/(t+x)$ belongs to $L^2((0, \infty ); dt)$
with $\Vert \phi_t\Vert_{L^2}^2=\pi /(2\sqrt x)$ and
satisfies 
$$\eqalignno{{{\sqrt{x}-\sqrt{y}}\over{x-y}}&={{1}\over{\pi}}\int_0^\infty {{ \sqrt{t}dt}\over
{(x+t)(y+t)}}\cr
&={{1}\over{\pi}}\int_0^\infty \phi_t(x)\phi_t(y)\, dt.&(5.13)}$$
\noindent We observe that $x^{1/4}u(x)\phi_t(x)$ belongs to $L^2((1, \infty ); dx/x)$ 
and since 
$$\eqalignno{\Bigl(x{{\partial}\over{\partial x}}+y{{\partial}\over{\partial
y}}\Bigr)W(x,y)&=-{{1}\over{4}}
(xy)^{1/4}u(x)u(y)\cr
&\quad+{{\kappa}\over{2\pi}}(xy)^{1/4}u(x)u(y)\int_0^\infty 
\phi_t(x)\phi_t(y)\, dt&(5.14)\cr}$$
\noindent we have
$$\eqalignno{W(x,y)&={{1}\over{4}}\int_1^\infty
u(xs)u(ys)(xy)^{1/4}{{ds}\over{s^{1/2}}}\cr
&\quad -{{\kappa}\over{2\pi}}\int_1^\infty \int_0^\infty \phi_t(sx)\phi_t(sy)
u(sx)u(sy)(xy)^{1/4}{{dtds}\over{s^{1/2}}}.&(5.15)}$$
\noindent Hence $W=\Gamma_\Phi^*\Gamma_\Phi$, where $K={\bf R}\oplus L^2((0, \infty );
dt)$ and $\Phi :(0, \infty )\rightarrow K$ is
$$\Phi (x)=2^{-1}x^{1/4}u(x)\oplus (-\kappa /2\pi
)^{1/2}x^{1/4}u(x)\phi_t(x).\eqno(5.16)$$
\rightline{$\square$}\par
\vskip.05in
\noindent {\bf Remarks 5.4.} (i) The condition $\kappa \leq 0$ in Proposition
5.3 excludes the case of the
associated Laguerre functions $w(x)=x^{(1+\alpha )/2}e^{-x/2}L_n^{(\alpha )}(x)$, 
where $L^{(\alpha )}_n$ with $(n=0,1, 2, \dots )$ are the associated 
Laguerre polynomials as in [14] and 
[7, 9.237].\par 
\noindent (ii) The Laplace transforms of $\sqrt {x} K_\nu (2\sqrt {x})$ and 
$\sqrt {x} J_\nu (2\sqrt {x})$ may be expressed in terms of Whittaker's functions.\par  
\noindent (iii) Bessel's equation may be transformed into the typical confluent
hypergeometric equation as in [8, p. 228]. \par 
\vskip.05in
\noindent {\bf Acknowledgement.} I am grateful to the referee for
pointing out reference [20].\par
\vskip.05in
\noindent {\bf References}\par
\vskip.05in
\noindent [1] G. Aubrun, A sharp small deviation inequality for the largest eigenvalue
of a random matrix, {Springer Lecture Notes 
 Lecture Notes in Math.,} {1857}, Springer, Berlin, 2005.\par 
 \noindent [2] G. Blower, {Operators associated with the soft and hard edges from unitary
ensembles}, {J. Math. Anal. Appl.} 2007. doi:10.1016/j.jmaa.2007.03.084.\par
\noindent [3] A. Borodin and G. Olshanski, 
Distributions on partitions, point processes, and the hypergeometric kernel, { 
Comm. Math. Phys.} {211} (2000), 335--358.\par
\noindent [4] A. Borodin and G. Olshanski, Asymptotics of Plancherel type random
partitions, {\sl arXiv:math.PR/0610240v2.}\par
\noindent [5] P.A. Deift, A.R. Its and X. Zhou, 
A Riemann--Hilbert approach to asymptotic problems arising in the theory of 
random matrix models, and also in the theory of integrable statistical mechanics,  
{Ann. of Math.} (2)  {146}  (1997), 149--235.\par
\noindent [6] F.G. Gantmacher, The Theory of Matrices, Volume II, Chelsea, 
New York, 1959.\par
\noindent [7] I.S. Gradshteyn and I.M. Ryzhik, {Table of Integrals, Series and
Products},
Fifth Edition, Academic Press, 1994.\par
\noindent [8] E. Hille, Ordinary Differential Equations in the Complex Domain, Wiley, London, 1976.\par 
\noindent [9] R.A. Horn and C.R. Johnson, Topics in Matrix Analysis, Cambridge University Press,
1991.\par
\noindent [10] E.L. Ince, {Ordinary Differential Equations},  
Dover Publications, London, 1956.\par
\noindent [11] A. McCafferty, Integrable operators and squares of 
Hankel matrices,\par
\noindent {\sl arXiv:math/0707.1474v2.}\par
\noindent [12] A.N. Megretski\u\i, V.V. Peller and S.R. Treil, 
The inverse spectral problem for self-adjoint Hankel operators, {Acta Math.} 
{174} (1995), 241--309.\par
\noindent [13] V. Peller, {Hankel Operators and Their
Applications}, Springer, New York, 2003.\par
\noindent [14] G. Sansone, {Orthogonal functions,} Interscience, New York, 1959.\par
\noindent [15] A.G. Soshnikov, Determinantal random point fields, 2000, 
{\sl arXiv.org:math/0002099.}\par
\noindent [16] G. Szeg\"o, {Orthogonal Polynomials,}
American Mathematical Society, New York, 1959.\par
\noindent [17] C.A. Tracy and H. Widom, {Level spacing distributions and the Airy
kernel,} {Comm. Math. Phys.} {159} (1994), 
151--174.\par
\noindent [18] C.A. Tracy and H. Widom, {Level spacing distributions and the Bessel
kernel,} {Comm. Math. Phys.} {161} (1994), 
289--309.\par
\noindent [19] C.A. Tracy and H. Widom, {Fredholm determinants, differential
equations and matrix models}, {Comm. Math. Phys.} {163} (1994), 33--72.\par
\noindent [20] C.A. Tracy and H. Widom, {Systems of partial
differential equations for a class of operator determinants}, pp.
381--388 in
{Partial Differential Operators and Mathematical Physics: Advances
and Applications, vol. 78}, Birkh\"auser Verlag, Berlin, 1995.\par  
\noindent [21] H.L. Turrettin, {Reduction of ordinary differential equations to the Birkhoff
canonical form,} {Trans. Amer. Math. Soc.} {107} (1963), 485--507.\par
\vfill
\eject
\end